\documentclass[leqno,11pt]{amsart}
\usepackage{amsmath,amstext,amssymb,amsopn,amsthm,mathrsfs}
\usepackage{psfrag}
\usepackage[dvips]{color}
\usepackage{multicol}
\usepackage{picinpar}
\usepackage{enumerate}
\usepackage[hyphens]{url}        
\usepackage[breaklinks,colorlinks=true,linkcolor=blue,citecolor=blue, urlcolor=blue]{hyperref}                                                 
\usepackage[utf8x]{inputenc}
\usepackage{txfonts}
\usepackage{graphicx}
\usepackage[T1]{fontenc}
\usepackage{longtable}
\theoremstyle{plain}
\newtheorem{teo}{Theorem}[section]
\newtheorem{defi}{Definition}[section]
\newtheorem{corollary}{Corollary}[section]
\newtheorem{lemma}{Lemma}[section]

\newtheorem{obs}{Observation}[section]
\numberwithin{equation}{section}

\begin{document}

\title[Boundedness of General Alternative Gaussian Singular Integrals ] {Boundedness of General Alternative Gaussian Singular Integrals in Gaussian variable Lebesgue spaces}
\author{Eduard Navas}
\address{Departamento de Matem\'aticas, Universidad Nacional Experimental Francisco de Miranda, Punto Fijo, Venezuela.}
\email{enavas@correo.unefm.edu.ve}
   \author{Ebner Pineda}
\address{Departamento de Matem\'{a}tica,  Facultad de Ciencias Naturales y Matem\'aticas, ESPOL Guayaquil 09-01-5863, Ecuador.}
\email{epineda@espol.edu.ec}
\author{Wilfredo~O.~Urbina}
\address{Department of Mathematics and Actuarial Sciences, Roosevelt University, Chicago, IL,
   60605, USA.}
\email{wurbinaromero@roosevelt.edu}


\begin{abstract}
In a previous paper \cite{NavPinUrb}, we introduced a new class of Gaussian singular integrals, that we called the general alternative Gaussian singular integrals and study the boundedness of them on  $L^p(\gamma_d)$, $ 1 < p < \infty.$ In this paper, we study the boundedness of those operators on  Gaussian variable Lebesgue spaces under a certain additional condition of regularity on $p(\cdot)$ following  \cite{DalSco}.
\end{abstract}

\maketitle

\section{Introduction and Preliminaries}

The  {\em general Gaussian singular integrals}, generalizing the Gaussian higher order Riesz  transforms were initially introduced  by W. Urbina in \cite{ur1} and later, S. P\'erez \cite{pe} extend to  a much larger class.

\begin{defi} Given a $C^1$-function $F,$ satisfying the orthogonality condition
 \begin{equation} \label{Fort}\index{general Gaussian singular integrals}
 \int_{\mathbb{R}^d} F(x) \gamma_d(dx) = 0,
\end{equation}
and such that for every $\varepsilon >0,$ there exist constants,  $C_{\varepsilon}$ and $C_{\varepsilon}'$
such that 
\begin{equation} \label{condF}
|F(x)| \leq C_{\epsilon}e^{\epsilon|x|^2} \quad \mbox{and} \quad |\nabla F(x)| \leq C_{\epsilon}'
e^{\epsilon|x|^2}.
\end{equation}
Then, for each $m \in \Bbb N$ the generalized Gaussian singular integral  is defined as
\begin{equation}
T_{F,m} f(x) = \int_{{\mathbb{R}}^d} \int_0^1 \left(\frac{- \log r}{1-r^2}\right)^{\frac{m - 2}{2}}
r^{m} F\Big(\frac{y-rx}{\sqrt{1-r^2}}\Big) \frac{e^{- \frac{|y-rx|^2}{1-r^2}}}{(1-r^2)^{d/2+1}} \frac{dr}{r} f(y) dy. 
\end{equation}
$T_{F,m}$ can be written as
$$ T_{F,m} f(x) = \int_{{\mathbb{R}}^d}  {\mathcal K}_{F,m} (x,y)  f(y) dy,$$
denoting,
\begin{eqnarray}\label{kernelF}
\nonumber {\mathcal K}_{F,m} (x,y) &=&  \int_0^1 \left(\frac{- \log r}{1-r^2}\right)^{\frac{m - 2}{2}} r^{m-1} F\Big(\frac{y-rx}{\sqrt{1-r^2}}\Big) 
 \frac {e^{- \frac{|y-rx|^2}{1-r^2}}}{(1-r^2)^{d/2+1}} dr \\
 &=& \int_0^1 \varphi_m(r) F\Big(\frac{y-rx}{\sqrt{1-r^2}}\Big)  \frac{e^{- \frac{|y-rx|^2}{1-r^2}}}{(1-r^2)^{d/2+1}} dr \\
\nonumber & =&\frac12 \int_0^1 \psi_m(t) F\Big(\frac{y-\sqrt{1-t}\, x}{\sqrt{t}}\Big)  
\frac {e^{-u(t)}}{t^{d/2+1}} dt,\\
\nonumber
\end{eqnarray}
with $\varphi_m(r)= \left(\frac{- \log r}{1-r^2}\right)^{\frac{m - 2}{2}} r^{m-1};$ and taking the change of variables $t =1-r^2,$  with $\psi_m(t) = \varphi_m(\sqrt{1-t})/\sqrt{1-t},$  and $u(t) = \frac{|\sqrt{1-t} x -y|^2}{t}.$ 
\end{defi}

In \cite{pe}, S. P\'erez proved the boundedness of $T_{F,m}$ on Gaussian $L^p$ spaces,
\begin{teo}\label{strongGenGaussSingInt}
The operators $T_{F,m} $ are $L^p(\gamma_d)$ bounded for $1< p < \infty,$ that is to say there exists $C>0,$ depending only in $p$ and dimension such that 
\begin{equation}
\|T_{F,m} f\|_{p,\gamma} \leq C \| f\|_{p,\gamma},
\end{equation}
for any $f \in L^p(\gamma_d).$ \\
 \end{teo}
Regarding the weak $(1,1)$ boundeness with respect to the Gaussian variable, she proved a negative result,
\begin{teo}\label{NegRes} Let $\Omega_t = \Big\{ z \in \mathbb{R}^d: \min_{1\leq i \leq d} |z_i| \geq t\Big\}$ and $\varTheta(t) = \frac{\inf_{\Omega_t} F(z)}{t^2},$ if $\limsup_{t\to \infty} \varTheta(t) = \infty,$ then the operator $T_{F,m}$ is not of weak type $(1,1)$ with respect to the Gaussian measure.
\end{teo}

Also, she obtained a positive result that is contained in the following theorem. In order to get sufficient conditions on $F$ for the weak type $(1,1)$ of $T_{F,m},$ since it is known that  the Gaussian Riesz transform $\mathcal{R}_{\beta}$  for $|\beta|\geq 3$ are not weak $(1,1)$ with respect to the Gaussian measure. Thus, since the weak type is not true, the natural question is what  weights  can be put in order to get a weak type inequality. She got the the weight should be of the form $w(y) = 1+ |y|^{|\beta| -2}.$ 
Moreover, for every $0 < \epsilon <|\beta|-2,$ there exists a function $F \in L^1((1+|\cdot|^\epsilon) \gamma_d)$ such that $T_{F,m} f \notin L^{1,\infty} (\gamma_d),$ see \cite{forsc}. The weights $w$ that will be considered, in order to ensure that $T_{F,m}$ is bounded from $L^1(w\gamma_d)$ into $L^{1,\infty}(\gamma_d),$  depend on the function $\varPhi.$

\begin{teo}\label{PosRes} The operator $T_{F,m}$ maps continuously  $L^1(w\gamma_d)$ into $L^{1,\infty}(\gamma_d)$ with $w(y)= 1 \vee \max_{1 \leq t \leq|y|} \eta(t)$ and 
$$ \eta(t) = \begin{cases}
\varPhi(t)/t \quad \mbox{if}  & \;1 \leq m <2,\\
\varPhi(t)/t^2 \quad \mbox{if}  & \; m \geq 2,\\
\end{cases}
$$
\end{teo}

In a previous paper  \cite{NavPinUrb}, we introduced a new class of Gaussian singular integrals, the {\em general alternative Gaussian singular integrals} as follows:
\begin{defi}\label{AlterGaussSingInt} Given a $C^1$-function $F,$ satisfying the orthogonality condition
 \begin{equation} \label{Fort2}
  \int_{\mathbb{R}^d} F(x)\; \gamma_d(dx) = 0,
\end{equation}
and such that for every $\varepsilon >0,$ there exist constants,  $C_{\varepsilon}$ and $C_{\varepsilon}'$
such that 
\begin{equation} \label{condF2}
|F(x)| \leq C_{\epsilon}e^{\epsilon|x|^2} \quad \mbox{and} \quad |\nabla F(x)| \leq C_{\epsilon}'
e^{\epsilon|x|^2}.
\end{equation}
Then, for each $m \in {\mathbb N}$ the general alternative Gaussian singular integral is defined as
\begin{equation}
\overline{T}_{F,m} f(x) = \int_{{\mathbb{R}}^d} \int_0^1\left(\frac{- \log r}{1-r^2}\right)^{\frac{m - 2}{2}}
 r^{d-1}F\Big(\frac{x-ry}{\sqrt{1-r^2}}\Big) \frac{e^{- \frac{|y-rx|^2}{1-r^2}}}{(1-r^2)^{d/2+1}} dr f(y) dy. 
\end{equation}
Thus, $\overline{T}_{F,m}$ can be written as
$$\overline{T}_{F,m} f(x) = \int_{{\mathbb{R}}^d}  \overline{\mathcal{K}}_{F,m} (x,y)  f(y) dy,$$
where,
\begin{eqnarray}\label{kernelF}
\nonumber  \overline{\mathcal{K}}_{F,m} (x,y) &=&  \int_0^1 \left(\frac{- \log r}{1-r^2}\right)^{\frac{m - 2}{2}}  r^{d-1} F\Big(\frac{x-ry}{\sqrt{1-r^2}}\Big) 
 \frac {e^{- \frac{|y-rx|^2}{1-r^2}}}{(1-r^2)^{d/2+1}} dr \\
 &=& \int_0^1 \varphi_m(r)\; F\Big(\frac{x-ry}{\sqrt{1-r^2}}\Big)  \frac{e^{- \frac{|y-rx|^2}{1-r^2}}}{(1-r^2)^{d/2+1}} dr \\
\nonumber& =&\frac12 \int_0^1 \psi_m(t)\; F\Big(\frac{x-\sqrt{1-t}\, y}{\sqrt{t}}\Big)  
\frac {e^{-u(t)}}{t^{d/2+1}} dt,\\
\nonumber
\end{eqnarray}
with $\varphi_m(r)= \left(\frac{- \log r}{1-r^2}\right)^{\frac{m - 2}{2}} r^{d-1};$ and after taking the change of variables $t =1-r^2,$   $\psi_m(t) = \varphi_m(\sqrt{1-t})/\sqrt{1-t},$  and $u(t) = \frac{|y- \sqrt{1-t} x|^2}{t}.$ 
\end{defi}
Additionally, in \cite{NavPinUrb} it was proved the boundedness of them in  $L^p(\gamma_d)$, for $d >1$ and $ 1 < p < \infty,$

\begin{teo}\label{strongGenGaussSingInt3}
For $d>1$, the operators $ \overline{T}_{F,m} $ are $L^p(\gamma_d)$ bounded for $1< p < \infty,$ that is to say there exists $C>0,$ depending only in $p$ and dimension such that 
\begin{equation}
\| \overline{T}_{F,m} f\|_{p,\gamma} \leq C \| f\|_{p,\gamma},
\end{equation}
for any $f \in L^p(\gamma_d).$ \\
 \end{teo}
 
In \cite{aimarforzaniscot}, H. Aimar, L. Forzani and R. Scotto  obtained a surprising result: the alternative Riesz transforms  $\overline{\mathcal{R}}_\beta$  are weak type $(1,1)$  for all multi-index $\beta$, i. e.  independently of their orders which is a contrasting fact with respect to the anomalous behavior of the higher order Riesz transforms $\mathcal{R}_\beta.$ For the general alternative Gaussian singular integrals $ \overline{T}_{F,m}$ we also  proved

\begin{teo} \label{altrieszteo} For $d > 1,$ there exists a constant $ C
$ depending only on $d$ and $m$ such that for all $
\lambda >0 $ and $ f \in L^1(\gamma_d) $, we have
$$\gamma_d \Big(\Big\{ x\in \mathbb{R}^d:\; \overline{T}_{F,m}(x) > \lambda \Big\}\Big) \le \frac{C}{\lambda}
 \int_{\mathbb{R}^d} |f(y)|\gamma_d(dy).$$
\end{teo}

On the other hand, in \cite{DalSco} E. Dalmasso and R. Scotto proved the boundedness of the general Gaussian singular integrals $T_{F,m}, $ on  Gaussian variable Lebesgue spaces under certain condition of regularity on $p(\cdot).$ In order to understand their result we need to get more background on variable Lebesgue spaces with respect to a Borel measure in general and the Gaussian measure in particular.

As usual in what follows $C$ represents a constant that is not necessarily the same in each occurrence; also we will used the notation: given two functions $f$, $g$, the symbols $\lesssim$ and $\gtrsim$ denote, that there is a constant $c$ such that $f\leq cg$ and $cf\geq g$, respectively. When both inequalities are satisfied, that is, $f\lesssim g\lesssim f$, we will denote $f\approx g$. 

Any $\mu$-measurable function $p(\cdot):\mathbb{R}^{d}\rightarrow [1,\infty]$ is an exponent function; the set of all the exponent functions will be denoted by  $\mathcal{P}(\mathbb{R}^{d},\mu)$. For $E\subset\mathbb{R}^{d}$ we set $$p_{-}(E)=\text{ess}\inf_{x\in E}p(x) \;\text{and}\; p_{+}(E)=\text{ess}\sup_{x\in E}p(x).$$
We use the abbreviations $p_{+}=p_{+}(\mathbb{R}^{d})$ and $p_{-}=p_{-}(\mathbb{R}^{d})$.

\begin{defi}\label{deflogholder}
Let $E\subset \mathbb{R}^{d}$. We say that $\alpha(\cdot):E\rightarrow\mathbb{R}$ is locally log-H\"{o}lder continuous, and denote this by $\alpha(\cdot)\in LH_{0}(E)$, if there exists a constant $C_{1}>0$ such that
			\begin{eqnarray*}
				|\alpha(x)-\alpha(y)|&\leq&\frac{C_{1}}{log(e+\frac{1}{|x-y|})}
			\end{eqnarray*}
			for all $x,y\in E$. We say that $\alpha(\cdot)$ is log-H\"{o}lder continuous at infinity with base point at $x_{0}\in \mathbb{R}^{d}$, and denote this by $\alpha(\cdot)\in LH_{\infty}(E)$, if there exist  constants $\alpha_{\infty}\in\mathbb{R}$ and $C_{2}>0$ such that
			\begin{eqnarray*}
				|\alpha(x)-\alpha_{\infty}|&\leq&\frac{C_{2}}{log(e+|x-x_{0}|)}
			\end{eqnarray*}
			for all $x\in E$. We say that $\alpha(\cdot)$ is log-H\"{o}lder continuous, and denote this by $\alpha(\cdot)\in LH(E)$ if both conditions are satisfied.
			The maximum, $\max\{C_{1},C_{2}\}$ is called the log-H\"{o}lder constant of $\alpha(\cdot)$.
\end{defi}

\begin{defi}\label{defPdlog}
			We denote $p(\cdot)\in\mathcal{P}_{d}^{log}(\mathbb{R}^{d})$, if $\frac{1}{p(\cdot)}$ is log-H\"{o}lder continuous and  denote by $C_{log}(p)$ or $C_{log}$ the log-H\"{o}lder constant of $\frac{1}{p(\cdot)}$.
		\end{defi}

We will need the following technical result, for its proof see Lemma 3.26 in \cite{dcruz}.			
\begin{lemma}\label{lema3.26CU}
Let $\rho(\cdot):\mathbb{R}^{d}\rightarrow[0,\infty)$ be such that $\rho(\cdot)\in LH_{\infty}(\mathbb{R}^{d})$, $0<\rho_{\infty}<\infty$, and let $R(x)=(e+|x|)^{-N}$, $N>d/\rho_{-}$. Then there exists a constant $C$ depending on $d$, $N$ and the $LH_{\infty}$ constant of $r(\cdot)$ such that given any set $E$ and
any function $F$ with $0\leq F(y)\leq 1$ for $y\in E$,
\begin{eqnarray}
  \int_{E}F^{\rho(y)}(y)dy &\leq& C\int_{E}F(y)^{\rho_{\infty}}dy + \int_{E}R^{\rho_{-}}(y)dy,\label{3.26.1} \\
  \int_{E}F^{\rho_{\infty}}(y)dy &\leq& C\int_{E}F^{r(y)}(y)dy + \int_{E}R^{\rho^{-}}(y)dy.\label{3.26.2}
\end{eqnarray}
\end{lemma}

\begin{defi}
For a $\mu$-measurable function $f:\mathbb{R}^{d}\rightarrow \mathbb{R}$, we define the modular \begin{equation}
\rho_{p(\cdot),\mu}(f)=\displaystyle\int_{\mathbb{R}^{d}\setminus\Omega_{\infty}}|f(x)|^{p(x)}\mu(dx)+\|f\|_{L^{\infty}(\Omega_{\infty},\mu)},
\end{equation}
and the norm
\begin{equation}
\|f\|_{L^{p(\cdot)}(\mathbb{R}^{d},\mu)}=\inf\left\{\lambda>0:\rho_{p(\cdot),\mu}(f/\lambda)\leq 1\right\}.
\end{equation}

\end{defi}

\begin{defi} The variable exponent Lebesgue space on $\mathbb{R}^{d}$, $L^{p(\cdot)}(\mathbb{R}^{d},\mu)$ consists on those $\mu\_$measurable functions $f$ for which there exists $\lambda>0$ such that $\rho_{p(\cdot),\mu}\left(\frac{f}{\lambda}\right)<\infty,$ i.e.
\begin{equation*}
L^{p(\cdot)}(\mathbb{R}^{d},\mu) =\left\{f:\mathbb{R}^{d}\to \mathbb{R}: f \; \text{measurable } \; \rho_{p(\cdot),\mu}\left(\frac{f}{\lambda}\right)<\infty, \; \text{for some} \;\lambda>0\right\}.
\end{equation*}

\end{defi}
It is well known that, if $p(\cdot) \in L H\left(\mathbb{R}^{d}\right)$ with $1<p_{-} \leq p^{+}<\infty$ the classical Hardy-Littlewood maximal function $\mathcal{M}$ is bounded on the variable Lebesgue space $L^{p(\cdot)},$ see \cite{dcruz1}. However, it is known that even though these are the sharpest possible point-wise conditions, they are not necessary. In \cite{LibroDenHarjHas} a necessary and sufficient condition is given for the $L^{p(\cdot)}$-boundedness of $\mathcal{M},$ but it is not an easy to work condition. The class $L H(\mathbb{R}^{d})$ is also sufficient for the boundedness on $L^{p(\cdot)}$-spaces of classical singular integrals of Calderón-Zygmund type, see \cite[Theorem 5.39]{dcruz}.\\

If $\mathcal{B}$ is a family of balls (or cubes) in $\mathbb{R}^{d}$, we say that $\mathcal {B}$ is $N$-finite if it has bounded overlappings for $N$, this is $\displaystyle\sum_{B\in\mathcal{B}}\chi_{B}(x)\leq N$ for all $x\in\mathbb {R}^{d}$; in other words, there is only $N$ balls (resp cubes) that intersect at the same time.\\

The following definition was introduced for the first time by Berezhno\v{\i} in \cite{Berez}, defined for family of disjoint balls or cubes. In the context of variable spaces, it has been considered in \cite{LibroDenHarjHas}, allowing the family to have bounded overlappings.\\
\begin{defi}
Given an exponent $p(\cdot)\in\mathcal{P}(\mathbb{R}^{d})$, we will say that $p(\cdot)\in\mathcal{G}$, if for every family of balls (or cubes) $\mathcal{B}$ which is $N$-finite,
\begin{eqnarray*}
  \sum_{B\in\mathcal{B}}||f\chi_{B}||_{p(\cdot)}||g\chi_{B}||_{p^{'}(\cdot)} &\lesssim& ||f||_{p(\cdot)}||g||_{p^{'}(\cdot)}
\end{eqnarray*}
for all functions $f\in L^{p(\cdot)}(\mathbb{R}^{d})$ and $g\in L^{p^{'}(\cdot)}(\mathbb{R}^{d})$. The constant only depends on N.
\end{defi}

\begin{lemma}[Teorema 7.3.22 in \cite{LibroDenHarjHas}]\label{implication1}
If $p(\cdot)\in LH(\mathbb{R}^{d})$, then $p(\cdot)\in\mathcal{G}$
\end{lemma}

We will consider only variable Lebesgue  spaces with respect to the Gaussian measure $\gamma_d,$ $L^{p(\cdot)}(\mathbb{R}^{d},\gamma_d).$  The next condition was introduced by E. Dalmasso and R. Scotto in \cite{DalSco}.

\begin{defi}\label{defipgamma}
Let $p(\cdot)\in\mathcal{P}(\mathbb{R}^{d},\gamma_{d})$, we say that $p(\cdot)\in\mathcal{P}_{\gamma_{d}}^{\infty}(\mathbb{R}^{d})$ if there exist constants $C_{\gamma_{d}}>0$ and $p_{\infty}\geq1$ such that
\begin{equation}
   |p(x)-p_{\infty}|\leq\frac{C_{\gamma_{d}}}{|x|^{2}},
\end{equation}
for $x\in\mathbb{R}^{d}\setminus\{(0,0,\ldots,0)\}.$
\end{defi}

\begin{obs}\label{obs4.1}
If $p(\cdot)\in\mathcal{P}_{\gamma_{d}}^{\infty}(\mathbb{R}^{d})$, then $p(\cdot)\in LH_{\infty}(\mathbb{R}^{d})$
\end{obs}
\begin{lemma}\label{lemaequiPgamma}
If $1<p_{-}\leq p_{+}<\infty,$ the following statements are equivalent
\begin{itemize}
  \item [(i)] $p(\cdot)\in\mathcal{P}_{\gamma_{d}}^{\infty}(\mathbb{R}^{d})$
  \item [(ii)] There exists $p_{\infty}>1$ such that
  \begin{eqnarray}
    C_{1}^{-1}\leq e^{-|x|^{2}(p(x)/p_{\infty}-1)}\leq C_{1} &\;\;\hbox{and}\;\;& C_{2}^{-1}\leq e^{-|x|^{2}(p^{'}(x)/p^{'}_{\infty}-1)}\leq C_{2},
  \end{eqnarray}
  for all $x\in\mathbb{R}^{d}$, where $C_{1}=e^{C_{\gamma_{d}}/p_{\infty}}$ and $C_{2}=e^{C_{\gamma_{d}}(p_{-})^{'}/p_{\infty}}$.
\end{itemize}
\end{lemma}

 Definition \ref{defipgamma} with Observation \ref{obs4.1} and Lemma \ref{lemaequiPgamma} end up strengthening the regularity conditions on the exponent functions $p(\cdot)$ to obtain the boundedness of the Ornstein-Uhlenbeck semigroup $\{T_{t}\}$, see \cite{MorPiUrb}. As a consequence of Lemma \ref{implication1}, we have
\begin{corollary}\label{solapamientoacotadoG}
If $p(\cdot)\in\mathcal{P}_{\gamma_{d}}^{\infty}(\mathbb{R}^{d})\cap LH_{0}(\mathbb{R}^{d})$, then $p(\cdot)\in\mathcal{G}$
\end{corollary}

As it has been mentioned already, in \cite{DalSco} E. Dalmasso and R. Scotto proved the boundedness of \(T_{F, m}\) on Gaussian variable Lebesgue spaces under the additional condition of regularity $p(\cdot)\in\mathcal{P}_{\gamma_{d}}^{\infty}(\mathbb{R}^{d})$.
 \begin{teo}\label{strongGenGaussSingInt}
Let $p(\cdot)\in\mathcal{P}_{\gamma_{d}}^{\infty}(\mathbb{R}^{d})\cap LH_{0}(\mathbb{R}^{d})$ with $1<p_{-}\leq p_{+}<\infty$. Then there exists a constant $C>0,$ depending only in $p$ and dimension such that 
\begin{equation}
\|T_{F,m} f\|_{p(\cdot),\gamma} \leq C \| f\|_{p(\cdot),\gamma},
\end{equation}
for any $f \in L^{p(\cdot)}(\gamma_d).$ \\
 \end{teo}
The main result in this paper is the proof, following the arguments of  Dalmasso and Scotto \cite{DalSco}, that the general alternative Gaussian singular integrals \(\overline{T}_{F, m}\) are also bounded on Gaussian variable Lebesgue spaces under the same condition of regularity on \(p(\cdot)\) considered by Dalmasso and Scotto.

\begin{teo}\label{strongGenGaussSingInt4}
Let $d >1$ and $p(\cdot)\in\mathcal{P}_{\gamma_{d}}^{\infty}(\mathbb{R}^{d})\cap LH_{0}(\mathbb{R}^{d})$ with $1<p_{-}\leq p_{+}<\infty$. Then there exists a constant $C>0,$ depending only in $p$ and dimension such that 
\begin{equation}
\| \overline{T}_{F,m} f\|_{p(\cdot),\gamma} \leq C \| f\|_{p(\cdot),\gamma},
\end{equation}
for any $f \in L^{p(\cdot)}(\gamma_d).$ \\
 \end{teo}

\section{Proof of the main result.}

We are ready for the proof of our main result, Theorem \ref{strongGenGaussSingInt4}. As usual we split operator \(\overline{T}_{F, m}\) into a local and a global part,
$$
\begin{aligned}
\overline{T}_{F, m} f(x) &=C_{d} \int_{|x-y|<d m(x)} \overline{\mathcal{K}}_{F, m}(x, y) f(y) d y+C_{d} \int_{|x-y| \geq d m(x)} \overline{\mathcal{K}}_{F, m}(x, y) f(y) d y \\
&=\overline{T}_{F, m, L} f(x)+\overline{T}_{F, m, G,} f(x)
\end{aligned}
$$
where
$$
\overline{T}_{F, m, L} f(x)=\overline{T}_{F, m}\left(f \chi_{B_{h}(\cdot)}\right)(x)
$$
is the {\em local part} and
$$
\overline{T}_{F, m, G} f(x)=\overline{T}_{F, m}\left(f \chi_{B_{h}^{c}(\cdot)}\right)(x)
$$
is the {\em global part} of $\overline{T}_{F, m},$ and
$$B_h(x) = B(x, C_d m(x)) = \{y \in \mathbb{R}^d: |y-x| < C_d m(x) \},$$
with  $m(x) = 1\wedge \frac{1}{|x|},$ is an admissible (hyperbolic) ball for the Gaussian measure, see \cite[Chapter 1]{urbina2019}.
\begin{proof} We need to bound the local and the global part.
\begin{enumerate}
\item[i)] 
The study of the local part $\overline{T}_{F, m, L}$ is similar to the one done in the proof of \cite[Theorem 1.2]{NavPinUrb}, see also \cite[Lemma 3.2]{DalSco}, obtaining the inequality
\begin{eqnarray}\label{ineq1}
\nonumber \left|\overline{T}_{F, m, L} f(x)\right|&=&\left|\overline{T}_{F, m} f\left(\chi_{B_{h}(\cdot)}\right)(x)\right|=\left|\int_{B_{h}(x)} \overline{\mathcal{K}}_{F, m}(x-y) f(y) dy\right|\\
&\leq& \sum_{B \in \mathcal{F}} \left|T\left(f \chi_{\hat{B}(\cdot)}\right)(x)\right|+\mathcal{M}\left(f \chi_{\hat{B}(\cdot)}\right)(x).
\end{eqnarray}
 where \(\mathcal{M}(g)\) is the classical Hardy-Littlewood maximal function of the function \(g\), and
$$
T f(x)= \text{p.v.} \int_{\mathbb{R}^{n}} \mathcal{K}(x-y) f(y) dy
$$
is a (convolution type) Calder\'on-Zygmund operator with kernel
$$
\mathcal{K}(x)=\int_{0}^{\infty} F\left(-\frac{x}{t^{1 / 2}}\right) e^{-|x|^{2} / t} \frac{d t}{t^{d / 2+1}}
$$
and \(\mathcal{F}\) is a countable family of admissible balls such that satisfies
the conditions of (covering) Lemma 4.3 of \cite{urbina2019} in particular, \(\mathcal{F}\) verify
\begin{enumerate}
\item[i)] For each $B \in\mathcal{F}$ let  \(\tilde{B}=2 B,\) then, the family of those balls $\tilde{\mathcal{F}}= \{B(0,1),\{\tilde{B}\}_{B \in \mathcal{F}}\}$ is a covering of \(\mathbb{R}^{d}\);
\item[ii)] \(\mathcal{F}\) has a bounded overlaps property;
\item[iii)] Every ball \(B \in \mathcal{F}\) is contained in an admissible ball, and therefore for
any pair \(x, y \in B, e^{-|x|^{2}} \sim e^{-|y|^{2}}\) with constants independent of \(B\)
\item[iv)]  There exists a uniform positive constant \(C_{d}\) such that, if \(x \in B \in \mathcal{F}\) then \(B_h(x) \subset C_{d} B:=\hat{B} .\) Moreover, the collection \(\hat{\mathcal{F}}=\{\hat{B}\}_{B \in \mathcal{F}}\) also satisfies the properties ii) and iii).\\
\end{enumerate}

Then, for \(f \in L^{p(\cdot)}\left(\mathbb{R}^{n}, \gamma_d\right)\) we will use the norm on the dual space
\(L^{p^{\prime}(\cdot)}\left(\mathbb{R}^{n}, \gamma_d\right)\)

\begin{equation}\label{ineq2}
\left\|\overline{T}_{F, m}\left(f \chi_{B_h(\cdot)}\right)\right\|_{p(\cdot), \gamma} \leq 2 \sup _{\|g\|_{p^{\prime}(\cdot), \gamma} \leq 1} \int_{\mathbb{R}^{n}}\left|\overline{T}_{F, m}\left(f \chi_{B_h(\cdot)}\right)(x)\right||g(x)| \gamma_d(dx)
\end{equation}

Using the pointwise inequality (\ref{ineq1}) we split the integral as
\begin{eqnarray*}
\int_{\mathbb{R}^{d}}\left|\overline{T}_{F, m}\left(f \chi_{B_h(\cdot)}\right)(x)\right||g(x)| \gamma_d(dx) &\lesssim& 
\sum_{B \in \mathcal{F}} \int_{B}\left|T\left(f \chi_{\hat{B}(\cdot)}\right)(x)\right||g(x)| e^{-|x|^{2}} dx \\
&&\quad+\sum_{B \in \mathcal{F}} \int_{B} \mathcal{M}\left(f \chi_{\hat{B}(\cdot)}\right)(x)|g(x)| e^{-|x|^{2}} d x \\
&& \approx \sum_{B \in \mathcal{F}} e^{-\left|c_{B}\right|^{2}} \int_{B}\left|T\left(f \chi_{\hat{B}(\cdot)}\right)(x)\right||g(x)| d x \\
&& \quad+\sum_{B \in \mathcal{F}} e^{-\left|c_{B}\right|^{2}} \int_{B} \mathcal{M}\left(f \chi_{\hat{B}(\cdot)}\right)(x)|g(x)| d x,
\end{eqnarray*}
where \(c_{B}\) is the center of \(B\) and \(\hat{B}\) and we have used property iii) above, i.e. that over each ball of the family \(\mathcal{F},\) the values of \(\gamma_d\) are all equivalent. Applying Hölder's inequality with \(p(\cdot)\) and \(p^{\prime}(\cdot)\) with respect of the Lebesgue measure in each integral and the boundedness of \(T\) and \(\mathcal{M}\) on $L^{p(\cdot)}\left(\mathbb{R}^{n}\right),$ we get
\begin{eqnarray}\label{ineq3}
\nonumber \int_{\mathbb{R}^{d}}\left|\overline{T}_{F, m}\left(f \chi_{B_h(\cdot)}\right)(x)\right||g(x)| \gamma_{d}(d x) &\lesssim & \sum_{B \in \mathcal{F}} e^{\left|c_{B}\right|^{2}}\left\|T\left(f \chi_{\hat{B}(\cdot)}\right) \chi_{B}\right\|_{p(\cdot)}\left\|g \chi_{B}\right\|_{p^{\prime}(\cdot)} \\
\nonumber &+&\sum_{B \in \mathcal{F}} e^{-\left|c_{B}\right|^{2}}\left\|\mathcal{M}\left(f \chi_{\hat{B}(\cdot)}\right) \chi_{B}\right\|_{p(\cdot)}\left\|g \chi_{B}\right\|_{p^{\prime}(\cdot)} \\
\nonumber &\lesssim &  \sum_{B \in \mathcal{F}} \mathrm{e}^{-\left|c_{B}\right|^{2}}\left\|f \chi_{\hat{B}}\right\|_{p(\cdot)}\left\|g \chi_{B}\right\|_{p^{\prime}(\cdot)} \\
&=& \sum_{B \in \mathcal{F}} \mathrm{e}^{-\left|c_{B}\right|^{2} / p_{\infty}}\left\|f \chi_{\hat{B}}\right\|_{p(\cdot)} e^{-\left|c_{B}\right|^{2} / p_{\infty}^{\prime}}\left\|g \chi_{\hat{B}}\right\|_{p^{\prime}()}
\end{eqnarray}
since \(p \in P_{\gamma}^{\infty}\left(\mathbb{R}^{n}\right)\) and \(p^{-}>1, p^{\prime} \in P_{\gamma}^{\infty}\left(\mathbb{R}^{n}\right) .\) Thus, from Lemma \(1.4,\) for
every \(x \in \mathbb{R}^{d}\)
\begin{equation}\label{equiv}
e^{-|x|^{2}\left(p(x) / p_{\infty}-1\right)} \leq C_{1} \text { and } e^{-|x|^{2}\left(p^{\prime}(x) / p_{\infty}^{\prime}-1\right)} \leq C_{2}.
\end{equation}

Moreover, since the values of the Gaussian measure \(\gamma_{d}\) are all equivalent on each ball \(\hat{B}\), we have
\begin{eqnarray*}
\int_{\hat{B}}\left(\frac{|f(y)|}{e^{\left|c_{B}\right|^{2} / p_{\infty}}\left\|f \chi_{\hat{B}}\right\|_{p(\cdot), \gamma}}\right)^{p(y)} d y &\lesssim & \int_{\hat{B}}\left(\frac{|f(y)|}{\left\|f \chi_{\hat{B}}\right\|_{p(\cdot), \gamma}}\right)^{p(y)} e^{-|y|^{2}\left(p(y) / p_{\infty}-1\right)} \gamma_{d}(dy) \\
&\lesssim &\int_{\hat{B}}\left(\frac{|f(y)|}{\left\|f \chi_{\hat{B}}\right\|_{p(\cdot), \gamma}}\right)^{p(y)} \gamma_{d}(dy) \lesssim  1
\end{eqnarray*}
which yields
$$
e^{-\left|c_{B}\right|^{2} / p_{\infty}}\left\|f \chi_{\hat{B}}\right\|_{p(\cdot)} \lesssim \left\|f \chi_{\hat{B}}\right\|_{p(\cdot), \gamma}.
$$
Similarly, by applying the second inequality of (\ref{equiv}) we get
$$
e^{-\left|c_{B}\right|^{2} / p_{\infty}^{\prime}}\left\|g \chi_{\hat{B}}\right\|_{p^{\prime}(\cdot)} \lesssim \left\|g \chi_{\hat{B}}\right\|_{p^{\prime}(\cdot), \gamma}.
$$

Replacing both estimates in (\ref{ineq3}) we obtain
\begin{eqnarray*}
\int_{\mathbb{R}^{d}}\left|\overline{T}_{F, m}\left(f \chi_{B_h(\cdot)}\right)(x)\right||g(x)| \gamma_{d}(d x) &  \lesssim  &\sum_{B \in \mathcal{F}}\left\|f \chi_{\hat{B}}\right\|_{p(\cdot), \gamma}\left\|g \chi_{\hat{B}}\right\|_{p^{\prime}(\cdot), \gamma} \\
&=& \sum_{B \in \mathcal{F}}\left\|f \chi_{\hat{B}} e^{-|\cdot|^{2} / p(\cdot) |}\right\|_{p(\cdot)}\left\|g \chi_{\hat{B}} e^{-|\cdot|^{2} / p^{\prime}(\cdot)}\right\|_{p^{\prime}(\cdot)}
\end{eqnarray*}
since the family of balls \(\hat{\mathcal{F}}\) has bounded overlaps, from Corollary 1.1 applied to \(f e^{-|\cdot|^{2} / p(\cdot)} \in L^{p(\cdot)}(\mathbb{R}^{d})\) and \(g e^{-|\cdot|^{2} / p^{\prime}(\cdot)} \in L^{p^{\prime}(\cdot)}(\mathbb{R}^{d}),\) it follows that
$$
\int_{\mathbb{R}^{d}}\left|\overline{T}_{F,m}\left(f \chi_{B_h(\cdot)}\right)(x)\right||g(x)| \gamma_{d}(d x) \leqslant\|f\|_{p(\cdot), \gamma}\|g\|_{p^{\prime}(\cdot), \gamma}
$$
Taking the supremum over all functions \(g\) with \(\|g\|_{p^{\prime}(\cdot), \gamma} \leq 1,\) from (\ref{ineq2}) we get finally
$$
\left\|\overline{T}_{F, m, L} f\right\|_{p(\cdot), \gamma}=\left\|\overline{T}_{F, m}\left(f \chi_{B_h(\cdot)}\right)\right\|_{p(\cdot), \gamma} \leq C\|f\|_{p(\cdot), \gamma}.\\
$$

\item[ii)] For the global part, to handle the kernel $\overline{\mathcal{K}},$ we need the following result
 \begin{lemma}\label{acotacionkernel}
  Let us consider the kernel \(\overline{\mathcal{K}}_{F, m}(x, y)\) in the global region, i.e. \(y \notin B_{h}(x)=\left\{y \in \mathbb{R}^{d}:\right.\) \(\left.|y-x|<C_{d} m(x)\right\}.\) If \(a=|x|^{2}+|y|^{2}\) and \(b=2\langle x, y\rangle,\) we have the following inequalities:
\begin{enumerate}
\item[i)] If \(b \leq 0,\) for each \(0<\epsilon<1,\) there exists \(C_{\epsilon}>0\) such that
$$
\left|\overline{\mathcal{K}}_{F, m}(x, y)\right| \leq C_{\epsilon} e^{-|y|^{2}+\epsilon|x|^{2}}
$$
\item[ii)] If \(b>0,\) for each \(0<\epsilon<\frac{1}{d},\) there exists \(C_{\epsilon}>0\) such that
$$
\left|\overline{\mathcal{K}}_{F, m}(x, y)\right| \leq C_{\epsilon} e^{\epsilon\left(|x|^{2}-|y|^{2}\right)} \frac{e^{-(1-\epsilon) u\left(t_{0}\right)}}{t_{0}^{d / 2}}
$$
where \(t_{0}=2 \frac{\sqrt{a^{2}-b^{2}}}{a+\sqrt{a^{2}-b^{2}}}\) and \(u_{0}=\frac{1}{2}\left(|y|^{2}-|x|^{2}+|x+y||x-y|\right)\).
\end{enumerate}
\end{lemma}

 Now, for the study of the global part \(\overline{T}_{F, m, G}\) let us consider the set  \[E_{x}=\{y: \langle x, y \rangle >0\},\] and consider two cases:
\begin{itemize}
\item Case $b=2\langle x,y\rangle\leq 0:$ let  $0<\epsilon <\frac{1}{p_{+}},$ then for $f\in
L^{p(\cdot)}(\mathbb{R}^{d},\gamma_{d})$ with
$\|f\|_{p(\cdot),\gamma_{d}}=1$ we have, by Lemma \ref{acotacionkernel} i)
\begin{eqnarray*}
\int_{\mathbb{R}^{d}}\left(\int_{B_h^{c}(\cdot)\cap
E_{x}^{c}}|\overline{{\mathcal K}}_{F,m}
(x,y)||f(y)|dy\right)^{p(x)}\gamma_{d}(dx)&\leq
&\int_{\mathbb{R}^{d}}\left(\int_{B_h^{c}(\cdot)\cap
E_{x}^{c}}e^{-|y|^{2}+\epsilon|x|^{2}}|f(y)|dy\right)^{p(x)}\gamma_{d}(dx)\\
&\leq
&\int_{\mathbb{R}^{d}}\left(\int_{\mathbb{R}^{d}}|f(y)|\gamma_{d}(dy)\right)^{p_{-}\frac{p(x)}{p_{-}}}e^{\epsilon
p(x)|x|^{2}-|x|^{2}}dx\\
&\leq
&\int_{\mathbb{R}^{d}}\left(\int_{\mathbb{R}^{d}}|f(y)|^{p_{-}}\gamma_{d}(dy)\right)^{\frac{p(x)}{p_{-}}}e^{(\epsilon
p_{+}-1)|x|^{2}}dx.
\end{eqnarray*}
Since, by hypothesis,
$\int_{\mathbb{R}^{d}}|f(x)|^{p(x)}\gamma_{d}(dx)\leq 1$ and therefore,
\begin{eqnarray*}
\int_{\mathbb{R}^{d}}|f(y)|^{p_{-}}\gamma_{d}(dy)&\leq
&\int_{|f|>1}|f(y)|^{p(x)}\gamma_{d}(dy)+\int_{|f|\leq
1}\gamma_{d}(dy)\\
&\leq & 1+C_{d}.
\end{eqnarray*}
Then,
\begin{eqnarray*}
\int_{\mathbb{R}^{d}}\left(\int_{B_h^{c}(\cdot)\cap
E_{x}^{c}}|\overline{{\mathcal K}}_{F,m}
(x,y)||f(y)|dy\right)^{p(x)}\gamma_{d}(dx)&\leq &
\int_{\mathbb{R}^{d}}(1+C_{d})^{\frac{p(x)}{p_{-}}}e^{(\epsilon
p_{+}-1)|x|^{2}}dx\\
&\leq &
(1+C_{d})^{\frac{p_{+}}{p_{-}}}\int_{\mathbb{R}^{d}}e^{(\epsilon
p_{+}-1)|x|^{2}}dx=C_{p,d}.
\end{eqnarray*}
In other words, 
$$\|\overline{T}_{F,m}(f\chi_{B_h^{c}(\cdot)\cap E_{x}^{c}})
\|_{p(\cdot),\gamma_{d}}\leq C_{p,d},$$ for all $f\in
L^{p(\cdot)}(\mathbb{R}^{d},\gamma_{d})$ with
$\|f\|_{p(\cdot),\gamma_{d}}=1$.\\

\item Case $b=2\langle x,y\rangle>0$: for $f\in
L^{p(\cdot)}(\mathbb{R}^{d},\gamma_{d})$ with
$\|f\|_{p(\cdot),\gamma_{d}}=1$ we have

\begin{eqnarray*}
(I)&=& \int_{\mathbb{R}^{d}}\left(\int_{B_h^{c}(\cdot)\cap
E_{x}}|\overline{{\mathcal K}}_{F,m}
(x,y)||f(y)|dy\right)^{p(x)}\gamma_{d}(dx)\hspace{3cm}\\
&\leq
&C\int_{\mathbb{R}^{d}}\left(\int_{B_h^{c}(\cdot)\cap
E_{x}}e^{\epsilon(|x|^{2}-|y|^{2})}\frac{e^{-(1-\epsilon)u(t_{0})}}{t_{0}^{d/2}}|f(y)|dy\right)^{p(x)}e^{-|x|^{2}}(dx)\\
&=&C\int_{\mathbb{R}^{d}}\left(\int_{B_h^{c}(\cdot)\cap
E_{x}}\frac{e^{-(1-\epsilon)u(t_{0})}}{t_{0}^{d/2}}e^{\frac{|y|^{2}}{p(y)}}e^{\frac{-|x|^{2}}{p(x)}}e^{-\epsilon(|y|^{2}-|x|^{2})}|f(y)|e^{-\frac{|y|^{2}}{p(y)}}dy\right)^{p(x)}(dx).\\
\end{eqnarray*}
Now, using the inequality  $||y|^{2}-|x|^{2}|\leq |x+y||x-y|,\,$
and that on the global region, $|x+y||x-y|> d,$ as $b>0$,we have
\begin{eqnarray*}
\frac{e^{-(1-\epsilon)u(t_{0})}}{t_{0}^{d/2}}e^{\frac{|y|^{2}}{p(y)}}e^{\frac{-|x|^{2}}{p(x)}}e^{-\epsilon(|y|^{2}-|x|^{2})}&\leq
&C\frac{e^{-(1-\epsilon)u(t_{0})+\frac{|y|^{2}-|x|^{2}}{p_{\infty}}-\epsilon(|y|^{2}-|x|^{2})}}{t_{0}^{d/2}}\\
&=&\frac{C}{t_{0}^{d/2}}e^{-\frac{(1-\epsilon)}{2}\left(|y|^{2}-|x|^{2}+|x+y||x-y|\right)+(\frac{1}{p_{\infty}}-\epsilon)(|y|^{2}-|x|^{2})}\\
&=&\frac{C}{t_{0}^{d/2}}e^{\left(-\frac{(1-\epsilon)}{2}+\frac{1}{p_{\infty}}-\epsilon\right)(|y|^{2}-|x|^{2})}e^{-\frac{(1-\epsilon)}{2}\left(|x+y||x-y|\right)}\\
&\leq &C|x+y|^{d}e^{-\alpha_{\infty}|x+y||x-y|}\\
\end{eqnarray*}
where,
$$\alpha_{\infty}=\frac{(1-\epsilon)}{2}-\left|\frac{1}{p_{\infty}}-\epsilon-\frac{(1-\epsilon)}{2}\right|>0,\hspace{0.3cm}\mbox{if}\hspace{0.3cm}\epsilon<\frac{1}{p_{\infty}}.$$

Thus, we take
$0<\epsilon<min\{\frac{1}{p_{\infty}},\frac{1}{p_{+}}\}.$ Now, let us consider the kernel
$$P(x,y):=|x+y|^{d}e^{-\alpha_{\infty}|x+y||x-y|},$$ which is integrable in each variable (since it is symmetric), with constant
independent of $x$ and $y.$ Then,
\begin{eqnarray*}
&&\int_{B_h^{c}(\cdot)\cap
E_{x}}\frac{e^{-(1-\epsilon)u(t_{0})}}{t_{0}^{d/2}}e^{\frac{|y|^{2}}{p(y)}}e^{\frac{-|x|^{2}}{p(x)}}e^{-\epsilon(|y|^{2}-|x|^{2})}|f(y)|e^{-\frac{|y|^{2}}{p(y)}}dy\leq C
\int_{B^{c}(x)}P(x,y)|f(y)|e^{-\frac{|y|^{2}}{p(y)}}dy
\end{eqnarray*}
Set $A_{x}=\left\{y:\;\frac{d}{|x|}<|y-x|<\frac{1}{2}\right\}$ and
$C_{x}=B^{c}(x,1/2)=\left\{y:\;|y-x|>\frac{1}{2}\right\}.$ \\

Therefore, $B_h^{c}(x)\subset A_{x}\cup C_{x}$. Define
$$
  J_{1}=\int_{A_{x}\cap E_{x}}P(x,y)f(y)e^{-|y|^{2}/p(y)}dy,$$
  and
  $$ J_{2}=\int_{C_{x}\cap E_{x}}P(x,y)f(y)e^{-|y|^{2}/p(y)}dy.$$
Let us estimate  $J_{1}$ first. Observe that, if $y\in A_{x}$, 
 $\frac{3}{4}|x|\leq|y|\leq\frac{5}{4}|x|$ and then $|x|\approx|y|$ hence $|x|\approx|x+y|$, and thus
  \begin{eqnarray}\label{ineqJ}
\nonumber    J_{1}
    &\lesssim&\int_{\frac{d}{|x|}<|x-y|}|x|^{d}e^{-\alpha_{\infty}|x||x-y|}f(y)e^{-|y|^{2}/p(y)}dy\\
     &\lesssim& \mathcal{M}(fe^{-|\cdot|^{2}/p(\cdot)})(x).
  \end{eqnarray}
It is known that the Hardy-Littlewood maximal function is weak $(1,1)$ in in Lebesgue variable spaces, see \cite{dcruz1}, see also \cite{MorPiUrb}, then from the hypothesis on $p(\cdot)$ we get
$$
\|\mathcal{M}(fe^{-|\cdot|^{2}/p(\cdot)})\|_{p(\cdot)}
\lesssim\|fe^{-|\cdot|^{2}/p(\cdot)}\|_{p(\cdot)}=\|f\|_{p(\cdot),\gamma_{d}}=1,
$$
 and then
 $$\rho_{p(\cdot)}\left(\mathcal{M}(fe^{-|\cdot|^{2}
  /p(\cdot)})\right)\lesssim 1.$$
Now, in order to estimate $J_{2}$, we have
 $$ J_{2}\leq\|P(x,\cdot)\chi_{C_{x}}\|_{p^{'}(\cdot)}\leq C,$$
 for details see \cite{DalSco}. This implies that there exists a constant independent on $x$ such that,
$$
  J_{2} = \int_{C_{x}\cap E_{x}}P(x,y)f(y)e^{-|y|^{2}/p(y)}dy\leq C$$
  thus $$\frac{1}{C}\int_{C_{x}\cap E_{x}}P(x,y)f(y)e^{-|y|^{2}/p(y)}dy \leq 1.$$
We set $g(y)=f(y)e^{-|y|^{2}/p(y)}=g_{1}(y)+g_{2}(y)$, where
$g_{1}=g\chi_{\{g\geq1\}}$ and $g_{2}=g\chi_{\{g<1\}};$ applying
(\ref{ineqJ}), we have
\begin{eqnarray*}
  (I) &\lesssim& \int_{\mathbb{R}^{d}}\left(\int_{B^{c}(x)\cap E_{x}}P(x,y)f(y)e^{-|y|^{2}/p(y)}dy\right)^{p(x)}dx \\
   &\lesssim&\int_{\mathbb{R}^{d}}\left(J_{1}\right)^{p(x)}dx+\int_{\mathbb{R}^{d}}\left(J_{2}\right)^{p(x)}dx \\
   &\lesssim&\rho_{p(\cdot)}\left(M_{H-L}(fe^{-|\cdot|^{2}/p(\cdot)})\right) + \int_{\mathbb{R}^{d}}\left(\frac{1}{C}\int_{C_{x}\cap E_{x}}P(x,y)g_{1}(y)dy\right)^{p(x)}dx\\
   &&+\int_{\mathbb{R}^{d}}\left(\frac{1}{C}\int_{C_{x}\cap E_{x}}P(x,y)g_{2}(y)dy\right)^{p(x)}dx \\
  &\lesssim&1+(II)+(III)
\end{eqnarray*}
Now, we study the terms $(II)$ and $(III)$.
\begin{eqnarray*}
  (II) &=&  \int_{\mathbb{R}^{d}}\left(\frac{1}{C}\int_{C_{x}\cap E_{x}}P(x,y)g_{1}(y)dy\right)^{p(x)}dx\leq\int_{\mathbb{R}^{d}}\left(\int_{C_{x}\cap E_{x}}P(x,y)g_{1}(y)dy\right)^{p_{-}}dx
\end{eqnarray*}
On the other hand, using Lemma \ref{lema3.26CU} with
$G(x)=\frac{1}{C}\int_{C_{x}\cap E_{x}}P(x,y)g_{2}(y)dy\leq1$ and
applying the inequality (\ref{3.26.1}), we obtain
\begin{eqnarray*}
 (III)&=& \int_{\mathbb{R}^{d}}\left(\int_{C_{x}\cap E_{x}}\frac{1}{C}P(x,y)g_{2}(y)dy\right)^{p(x)}dx =\int_{\mathbb{R}^{d}}(G(x))^{p(x)}dx \\
  &\lesssim& \int_{\mathbb{R}^{d}}(G(x))^{p_{\infty}}dx + \int_{\mathbb{R}^{d}}\frac{dx}{(e+|x|)^{dp_{-}}}\\
  &=&\int_{\mathbb{R}^{d}}\left(\int_{C_{x}\cap E_{x}}P(x,y)g_{2}(y)dy\right)^{p_{\infty}}+C_{d,p}.
\end{eqnarray*}
Therefore
\begin{eqnarray*}
  (I)&\lesssim& \int_{\mathbb{R}^{d}}\left(\int_{C_{x}\cap E_{x}}P(x,y)g_{1}(y)dy\right)^{p_{-}}dx+\int_{\mathbb{R}^{d}}\left(\int_{C_{x}\cap E_{x}}P(x,y)g_{2}(y)dy\right)^{p_{\infty}}dx+C_{d,p}
\end{eqnarray*}
Now, in order to estimate the last two integrals, we apply H\"{o}lder's
inequality.
\begin{eqnarray*}
  &&\int_{\mathbb{R}^{d}}\left(\int_{C_{x}\cap E_{x}}P(x,y)g_{1}(y)dy\right)^{p_{-}}dx\leq\int_{\mathbb{R}^{d}}\left(\int_{\mathbb{R}^{d}}P(x,y)^{\frac{1}{p_{-}^{'}}}P(x,y)^{\frac{1}{p_{-}}}g_{1}(y)dy\right)^{p_{-}}dx\\
&\leq&\int_{\mathbb{R}^{d}}\left(\int_{\mathbb{R}^{d}}(P(x,y))^{p_{-}^{'}/p_{-}^{'}}dy\right)^{p_{-}/p_{-}^{'}}\left(\int_{\mathbb{R}^{d}}(P(x,y))^{p_{-}/p_{-}}g_{1}^{p_{-}}(y)dy\right)^{p_{-}/p_{-}}dx\\
&=&\int_{\mathbb{R}^{d}}\left(\int_{\mathbb{R}^{d}}P(x,y)dy\right)^{p_{-}/p_{-}^{'}}\left(\int_{\mathbb{R}^{d}}P(x,y)g_{1}^{p_{-}}(y)dy\right)dx\\
&\lesssim&\int_{\mathbb{R}^{d}}\int_{\mathbb{R}^{d}}P(x,y)g_{1}^{p_{-}}(y)dydx.
\end{eqnarray*}
Then, by Fubini's theorem we get,
\begin{eqnarray*}
  \int_{\mathbb{R}^{d}}\left(\int_{C_{x}\cap E_{x}}P(x,y)g_{1}(y)dy\right)^{p_{-}}dx &\lesssim&\int_{\mathbb{R}^{d}}\int_{\mathbb{R}^{d}}P(x,y)g_{1}^{p_{-}}(y)dydx\\
  &=&\int_{\mathbb{R}^{d}}g_{1}^{p_{-}}(y)\left(\int_{\mathbb{R}^{d}}P(x,y)dx\right)dy\\
  &\lesssim& \int_{\mathbb{R}^{d}}\left(g_{1}(y)\right)^{p(y)}dy \\
  &\lesssim& \int_{\mathbb{R}^{d}}f(y)^{p(y)}e^{-|y|^{2}}dy\lesssim\rho_{p(\cdot),\gamma_{d}}(f)
\end{eqnarray*}
To estimate the integral
$\int_{\mathbb{R}^{d}}\left(\int_{C_{x}\cap
E_{x}}P(x,y)g_{2}(y)dy\right)^{p_{\infty}}dx$, we proceed in
analogous way, but applying the H\"{o}lder's inequality to the
exponent $p_{\infty}$, and applying the inequality (\ref{3.26.2})
in Lemma \ref{lema3.26CU}. In consequence we obtain
\begin{eqnarray*}
  \int_{\mathbb{R}^{d}}\left(\int_{C_{x}\cap E_{x}}P(x,y)g_{2}(y)dy\right)^{p_{\infty}}dx &\lesssim&\int_{\mathbb{R}^{d}}g_{2}^{p(y)}(y)dy +C_{d,p} \\
&\lesssim&\rho_{p(\cdot),\gamma_{d}}(f)+C.
\end{eqnarray*}
Therefore,
\begin{eqnarray*}
  (I) &\lesssim& \int_{\mathbb{R}^{d}}\left(\int_{C_{x}\cap E_{x}}P(x,y)g_{1}(y)dy\right)^{p_{-}}dx+\int_{\mathbb{R}^{d}}\left(\int_{C_{x}\cap E_{x}}P(x,y)g_{2}(y)dy\right)^{p_{\infty}}dx+C_{d,p}\\
  &\leq&2\rho_{p(\cdot),\gamma_{d}}(f)+C_{d,p}\\
\end{eqnarray*}
Hence, we obtain that 
$$\| \overline{T}_{F,m} f(\chi_{B_h^{c}(\cdot)\cap
E_{x}})\|_{p(\cdot),\gamma}\leq C,$$
for
$\|f\|_{p(\cdot),\gamma_{d}}=1.$\\

Putting together both cases we get
$$\| \overline{T}_{F,m} f(\chi_{B_h^{c}(\cdot)})\|_{p(\cdot),\gamma}\leq C,$$
for
$\|f\|_{p(\cdot),\gamma_{d}}=1.$
 Then by homogeneity of the norm the result holds for all function $f\in L^{p(\cdot)}(\mathbb{R}^{d},\gamma_{d})$. Now the proof  of the theorem is complete.
\end{itemize}
\end{enumerate}
\end{proof}

\end{document}